\documentclass[12pt]{amsart}

\usepackage{amsmath}
\usepackage{amssymb}
\usepackage{array}

\newtheorem{theorem}{Theorem}%[section]

\theoremstyle{definition}

\begin{document}

\title{Simultaneous generation for zeta values by the Markov-WZ method}

% author one information
\author{Kh.~Hessami Pilehrood$^{1}$}

\address{Institute for Studies in Theoretical Physics and Mathematics
(IPM), Tehran, Iran} \curraddr{Mathemetics Department, Faculty of
Science, Shahrekord University, Shahrekord, P.O. Box 115, Iran.}
\email{hessamik@ipm.ir, hessamit@ipm.ir, hessamit@gmail.com}
\thanks{$^1$ This research was in part supported by a grant
from IPM (No. 86110025)}

% author two information
\author{T.~Hessami Pilehrood$^2$}
%\address{Institute for Studies in Theoretical Physics and Mathematics
%(IPM), Tehran, Iran}
%\curraddr{
%Mathemetics Department, Shahrekord University, Shahrekord,
%P.O. Box 115, Iran.
%}
%\email{hessamit@gmail.com}
\thanks{$^2$ This research was in part supported by a grant
from IPM (No. 86110020)}

%Primary. Secondary.}

\date{}


\begin{abstract}
By application of the Markov-WZ method, we prove a more general
form of a bivariate generating function identity containing, as
particular cases, Koecher's and Almkvist-Granville's Ap\'ery-like
formulae for odd zeta values. As a consequence, we get a new
identity producing Ap\'ery-like series for all $\zeta(2n+4m+3),$
$n,m\ge 0,$ convergent at the geometric rate with ratio $2^{-10}.$
\end{abstract}

\maketitle

\section{Introduction}
\label{intro}

The Riemann zeta function is defined by the series
$$
\zeta(s)=\sum_{n=1}^{\infty}\frac{1}{n^s}, \qquad \mbox{for} \quad
{\rm Re}(s)>1.
$$
Ap\'ery's irrationality proof of $\zeta(3)$  \cite{po} operates
with the faster convergent series
$$
\zeta(3)=\frac{5}{2}\sum_{k=1}^{\infty}\frac{(-1)^{k-1}}{k^3\binom{2k}{k}}
$$
first obtained by A.~A.~Markov in 1890 \cite{ma}. The general
formula giving analogous series for all $\zeta(2s+3),$ $s\ge 0,$
was proved by Koecher \cite{ko} (and independently in an expanded
form by  Leshchiner \cite{le})
\begin{equation}
\sum_{s=0}^{\infty}\zeta(2s+3)x^{2s}=\sum_{k=1}^{\infty}\frac{1}{k(k^2-x^2)}=
\frac{1}{2}\sum_{k=1}^{\infty}\frac{(-1)^{k-1}}{k^3\binom{2k}{k}}\,\,\frac{5k^2-x^2}%
{k^2-x^2}\,\prod_{m=1}^{k-1}\left(1-\frac{x^2}{m^2}\right).
\label{eq01}
\end{equation}
The similar identity generating fast convergent series for all
$\zeta(4s+3),$ $s\ge 0,$ which for $s>1$ are different from
Koecher's result (\ref{eq01})  was experimentally discovered in
\cite{bb} and proved by G.~Almkvist and A.~Granville in
\cite{algr}
\begin{equation}
\sum_{s=0}^{\infty}\zeta(4s+3)x^{4s}=\sum_{k=1}^{\infty}\frac{k}%
{k^4-x^4}=\frac{5}{2}\sum_{k=1}^{\infty}\frac{(-1)^{k-1}}%
{\binom{2k}{k}}\frac{k}{k^4-x^4}\prod_{m=1}^{k-1}\left(\frac{m^4+4x^4}%
{m^4-x^4}\right). \label{eq02}
\end{equation}
There exists a bivariate unifying formula for identities
(\ref{eq01}) and (\ref{eq02})
\begin{equation}
\sum_{k=1}^{\infty}\frac{k}{k^4-x^2k^2-y^4}=\frac{1}{2}\sum_{k=1}^{\infty}
\frac{(-1)^{k-1}}{k\binom{2k}{k}}\frac{5k^2-x^2}{k^4-x^2k^2-y^4}\prod_{m=1}^{k-1}
\left(\frac{(m^2-x^2)^2+4y^4}{m^4-x^2m^2-y^4}\right), \label{eq03}
\end{equation}
which was first conjectured by H.~Cohen and then proved by
D.~Bradley  \cite{br} and, independently, by T.~Rivoal \cite{ri}.
This identity implies (\ref{eq01}) if $y=0,$ and gives
(\ref{eq02}) if $x=0.$ The proof of (\ref{eq03})  relies on
Borwein $\&$ Bradley's method \cite {bb} and consists of reduction
of (\ref{eq03}) to a finite non-trivial combinatorial identity
which can be proved on the base of Almkvist $\&$ Granville's  work
\cite{algr}.

In this short note, we prove a more general form of (\ref{eq03})
by application of the Markov-WZ  method (see \cite{ma,moze, momo,
ks}). Let us notice that Koecher's identity (\ref{eq01}) and
similar Leschiner's and Bailey-Borwein-Bradley's identities
\cite{le, bbb} generating rapidly convergent series for even zeta
values $\zeta(2s+2)$ can be proved with the help of  the  WZ
method (see \cite{he} for more details).

\begin{theorem}  \label{t1}
Let $a, b$ be complex numbers, with $|a|<1,$ $|b|<1.$ Then for
arbitrary complex numbers $A_0, B_0, C_0$ we have
$$
\sum_{k=1}^{\infty}\frac{A_0+B_0k+C_0k^2}{(k^2-a^2)(k^2-b^2)}=
\sum_{n=1}^{\infty}\frac{d_n}{\prod_{m=1}^n(m^2-a^2)(m^2-b^2)},
$$
with
\begin{equation*}
\begin{split}
d_n&=\frac{(-1)^{n-1}B_0(n-1)!(5n^2-a^2-b^2)}{2^{n+1}}\prod_{m=1}^{n-1}
\left(\frac{(m^2-(a^2+b^2))^2-4a^2b^2}{2m+1}\right) \\ &+
\frac{20n+5}{2(5n^2-2a^2-2b^2)}L_n+
\frac{35n^5-35n^3(a^2+b^2)+4n(3a^4+3b^4-4a^2b^2)}{4(5n^2-2a^2-2b^2)}L_{n-1},
\end{split}
\end{equation*}
where $L_n$ is a solution of the second order difference equation
\begin{equation*}
\begin{split}
4&(4n+3)(4n+5)(5n^2-2a^2-2b^2)L_{n+1}+2(n+1)p(n)L_{n} \\ &-n(n+1)
(5(n+1)^2-2a^2-2b^2)q(n)L_{n-1}=0,\quad n=1,2,\ldots
\end{split}
\end{equation*}
with initial conditions $L_0=C_0,$
$$
L_1=\left(\frac{1}{3}-\frac{2}{15}(a^2+b^2)\right)A_0+
\left(\frac{1}{6}(a^2+b^2)-\frac{2}{15}(a^4+b^4
-4a^2b^2)-\frac{1}{30}\right)C_0,
$$
and
\begin{equation}
\begin{split}
p&(n)=30n^7+105n^6+n^5(145-52(a^2+b^2))+n^4(100-130(a^2+b^2))
\\ &+n^3
(35-124(a^2+b^2)+56(a^4+b^4)-208a^2b^2)+n^2(5-56(a^2+b^2) \\
&+84(a^4+b^4) -312a^2b^2)
+n(80a^2b^2(a^2+b^2)-16(a^6+b^6)+48(a^4+b^4-3a^2b^2) \\
&-14(a^2+b^2))
+(10(a^2-b^2)^2-2(a^2+b^2)+40a^2b^2(a^2+b^2)-8(a^6+b^6)),
\label{p}
\end{split}
\end{equation}
\begin{equation}
\begin{split}
q(n)&=n^8-6n^6(a^2+b^2)+n^4(9(a^4+b^4)+30a^2b^2) \\
&-n^2(28a^2b^2(a^2+b^2)+4(a^6+b^6)) +16a^2b^2(a^2-b^2)^2.
\label{q}
\end{split}
\end{equation}
\end{theorem}

If in Theorem \ref{t1} we take $B_0=1,$  $A_0=C_0=0,$ then $L_n=0$
for all $n\ge 0$ and putting
\begin{equation}
 a^2=\frac{x^2+\sqrt{x^4+4y^4}}{2},
\qquad b^2=\frac{x^2-\sqrt{x^4+4y^4}}{2} \label{ab}
\end{equation}
 we get the
identity (\ref{eq03}).
%\begin{corollary}
%Let $a, b$ be complex numbers with $|a|<1,$ $|b|<1.$ Then
%\begin{equation}
%\sum_{k=1}^{\infty}\frac{k}{(k^2-a^2)(k^2-b^2)}=\frac{1}{2}
%\sum_{n=1}^{\infty}\frac{(-1)^{n+1}(5n^2-a^2-b^2)}%
%{\binom{2n}{n}n(n^2-a^2)(n^2-b^2)}\prod_{m=1}^{n-1}\left(
%\frac{(m^2-(a^2+b^2))^2-4a^2b^2}{(m^2-a^2)(m^2-b^2)}\right).
%\label{c}
%\end{equation}
%\end{corollary}

If $A_0=1,$ $B_0=C_0=a=b=0,$ we get the following series for
$\zeta(4)$  mentioned by Markov in \cite[p.18]{ma}:
$$
\zeta(4)=\sum_{n=1}^{\infty}\frac{1}{n!^4}\left(\frac{4n+1}{2n^2}
L_n+\frac{7n^3}{4}L_{n-1}\right),
$$
where $L_0=0,$ $L_1=1/3,$ and
$$
4(4n+3)(4n+5)L_{n+1}+2(n+1)^3(6n^3+9n^2+5n+1)L_n-n^7(n+1)^3L_{n-1}=0,
\quad n\ge 1.
$$

\begin{theorem} \label{t2}
Let $x, y$ be complex numbers such that $|x|^2+|y|^4<1.$ Then
\begin{equation}
\sum_{k=1}^{\infty}\frac{k}{k^4-x^2k^2-y^4}=\frac{1}{2}
\sum_{n=1}^{\infty}\frac{(-1)^{n-1}r(n)}{n\binom{2n}{n}}\,\,
\frac{\prod_{m=1}^{n-1}((m^2-x^2)^2+4y^4)}%
{\prod_{m=n}^{2n}(m^4-x^2m^2-y^4)}, \label{idt2}
\end{equation}
where
$$
r(n)=205n^6-160n^5+(32-62x^2)n^4+40x^2n^3+(x^4-8x^2-25y^4)n^2+10y^4n+y^4(x^2-2).
$$
\end{theorem}
Since
$$
\sum_{k=1}^{\infty}\frac{k}{k^4-x^2k^2-y^4}=\sum_{n=0}^{\infty}
\sum_{m=0}^{\infty}\binom{n+m}{n}\zeta(2n+4m+3)x^{2n}y^{4m},
$$
the formula (\ref{idt2}) generates Ap\'ery-like series for all
$\zeta(2n+4m+3),$ $n,m \ge 0,$ convergent at the geometric rate
with ratio $2^{-10}.$ So, for example, if $x=y=0,$ we get
Amdeberhan and Zeilberger's series \cite{az} for $\zeta(3)$
$$
\zeta(3)=\frac{1}{2}\sum_{n=1}^{\infty}\frac{(-1)^{n-1}(205n^2-160n+32)}%
{n^5\binom{2n}{n}^5}.
$$
If $y=0,$ we recover Theorem 4 from \cite{he}. If $x=0,$ we find,
in particular, the following expression for $\zeta(7):$
\begin{equation*}
\begin{split}
\zeta(7)&=\frac{1}{2}\sum_{n=1}^{\infty}\frac{(-1)^n(25n^2-10n+2)}{n^9\binom{2n}{n}^5}
\\[2pt]
&-\frac{1}{2}\sum_{n=1}^{\infty}\frac{(-1)^n(205n^2-160n+32)}{n^5\binom{2n}{n}^5}
\left(\sum_{m=1}^{2n}\frac{1}{m^4}+\sum_{m=1}^{n-1}\frac{3}{m^4}\right).
\end{split}
\end{equation*}

\newpage

\section{Proof of Theorem \ref{t1}.}

As usual,  let $(\lambda)_{\nu}$ be the Pochhammer symbol (or the
shifted factorial) defined by
\begin{equation*}
(\lambda)_{\nu}=\frac{\Gamma(\lambda+\nu)}{\Gamma(\lambda)}
=\begin{cases}
 1,     & \quad \nu=0; \\
\lambda(\lambda+1)\ldots (\lambda+\nu-1), & \quad \nu\in {\mathbb
N}.
\end{cases}
\end{equation*}
Let $a, b$ be complex numbers such that $|a|<1,$ $|b|<1.$ We start
with the kernel
$$
H(n,k)=\frac{(1+a)_k(1-a)_k(1+b)_k(1-b)_k}{(1+a)_{n+k+1}(1-a)_{n+k+1}
(1+b)_{n+k+1}(1-b)_{n+k+1}}
$$
and define two functions
\begin{equation*}
\begin{split}
F(n,k)&=H(n,k)(A(n)+B(n)(k+1)+C(n)(k+1)^2), \\
G(n,k)&=H(n,k)(D(n)+E(n)k+K(n)k^2+L(n)k^3),
\end{split}
\end{equation*}
with $7$ unknown coefficients $A(n), B(n), C(n), D(n), E(n), K(n),
L(n)$ as functions of $n.$ We require that
\begin{equation}
F(n+1,k)-F(n,k)=G(n,k+1)-G(n,k). \label{eq04}
\end{equation}
Since $F(n,k)$ and $G(n,k)$ are not proper hypergeometric, the
pair $(F,G)$ is a Markov WZ-pair (see \cite[p. 8]{moze}, \cite{z}
for definitions).

Substituting $F, G$ into (\ref{eq04}) and cancelling common
factors we get the following equation of degree $6$ in  a variable
$k:$
\begin{equation}
\begin{split}
&((n+k+2)^2-a^2)((n+k+2)^2-b^2)(A(n)+B(n)(k+1)+C(n)(k+1)^2) \\
&-A(n+1)-B(n+1)(k+1)- C(n+1)(k+1)^2=((n+k+2)^2-a^2) \\
&\times((n+k+2)^2-b^2)(D(n)+E(n)k+K(n)k^2+L(n)k^3)- ((k+1)^2-a^2)
\\ &\times((k+1)^2-b^2)(D(n)+E(n)(k+1)+K(n)(k+1)^2+L(n)(k+1)^3).
\label{eq05}
\end{split}
\end{equation}
To satisfy condition (\ref{eq04}), all the coefficients of the
powers of $(k+1)$ in the equation (\ref{eq05}) must vanish. This
leads to a system of first order linear recurrence equations with
polynomial coefficients for $A(n), B(n), C(n), D(n), E(n), K(n),
L(n)$
\begin{equation}
C(n)=(4n_1-3)L(n), \quad B(n)=(4n_1-2)K(n)-(10n_1^2-3)L(n), \quad
n_1=n+1, \label{eq06}
\end{equation}
\begin{equation}
A(n)=(4n_1-1)E(n)-(10n_1^2-1)K(n)+(20n_1^3+2n_1(a^2+b^2)-1)L(n),
\label{eq07}
\end{equation}
\begin{equation}
4D(n)=10n_1E(n)-(20n_1^2+2a^2+2b^2)K(n)+(35n_1^3+11n_1(a^2+b^2))L(n),
\label{eq08}
\end{equation}
\begin{equation}
\begin{split}
2(4n_1+1)L(n+1)&=2n_1(5n_1^2-2a^2-2b^2)E(n)-2n_1^2(15n_1^2-6(a^2+b^2))K(n)
\\ &+
n_1(63n_1^4-17n_1^2(a^2+b^2)-4(a^4+b^4))L(n), \label{eq09}
\end{split}
\end{equation}
\begin{equation}
\begin{split}
&2(4n_1+2)K(n+1)-2(10n_1^2+20n_1+7)L(n+1) \\
&=2n_1^2(5n_1^2-2(a^2+b^2))E(n)
-2(16n_1^5-8n_1^3(a^2+b^2)+n_1(a^2-b^2)^2)K(n) \\ &+
(70n_1^6-31n_1^4(a^2+b^2)+n_1^2(3a^4+3b^4-14a^2b^2))L(n),
\label{eq10}
\end{split}
\end{equation}
\begin{equation}
\begin{split}
& 4(20(n_1+1)^3+2(n_1+1)(a^2+b^2)-1)L(n+1) -4(10n_1^2+20n_1+9)K(n+1) \\
&+4(4n_1+3)E(n+1)
=(6n_1^5-6n_1^3(a^2+b^2)+16a^2b^2n_1)E(n) \\
&-(20n_1^6-22n_1^4(a^2+b^2)+2n_1^2(a^4+b^4 +22a^2b^2))K(n) \\
&+(45n_1^7-48n_1^5(a^2+b^2)+n_1^3(3a^4+3b^4+86a^2b^2)+8a^2b^2n_1(a^2+b^2))L(n).
\label{eq11}
\end{split}
\end{equation}
Multiplying equation (\ref{eq09}) by $n_1$ and subtracting from
(\ref{eq10}) we get
$$
2K(n+1)-7(n_1+1)L(n+1)=-\frac{n_1((n_1^2-
a^2-b^2)^2-4a^2b^2)}{2(2n_1+1)}(2K(n)-7n_1L(n)),
$$
which yields
$$
K(n)-\frac{7}{2}n_1L(n)=\frac{(-1)^n(2K(0)-7L(0))n!}{2^{n_1}}
\prod_{m=1}^n\left(\frac{(m^2-a^2-b^2)^2-4a^2b^2}{2m+1}\right).
$$
From (\ref{eq06}) it follows that $2K(0)=B(0)+7L(0)$ and
therefore,
\begin{equation}
K(n)=\frac{7}{2}n_1L(n)+\frac{(-1)^nB(0)n!}{2^{n_1}}\prod_{m=1}^n
\left(\frac{(m^2-a^2-b^2)^2-4a^2b^2}{2m+1}\right). \label{eq12}
\end{equation}
Substituting  (\ref{eq12}) into (\ref{eq09}) yields the formula
for $E(n)$
\begin{equation}
\begin{split}
 E(n)&=\frac{4n_1+1}{n_1(5n_1^2-2a^2-2b^2)}L(n+1)+
\frac{42n_1^4-25n_1^2(a^2+b^2)+4(a^4+b^4)}{2(5n_1^2-2a^2-2b^2)}L(n)
\\ &+\frac{3B(0)(-1)^nn_1!}{2^{n_1}}\prod_{m=1}^n\left(
\frac{(m^2-a^2-b^2)^2-4a^2b^2}{2m+1}\right). \label{eq13}
\end{split}
\end{equation}
Substituting  (\ref{eq12}) and (\ref{eq13}) into (\ref{eq08})
gives the formula for $D(n)$
\begin{equation}
\begin{split}
D(n)&=\frac{(40n_1+10)L(n+1)+
(35n_1^5-35n_1^3(a^2+b^2)+4n_1(3a^4+3b^4-4a^2b^2))L(n)}%
{4(5n_1^2-2a^2-2b^2)} \\
&+\frac{(-1)^nB(0)n!(5n_1^2-a^2-b^2)}{2^{n_1+1}}\prod_{m=1}^n
\left(\frac{(m^2-a^2-b^2)^2-4a^2b^2}{2m+1}\right). \label{eq14}
\end{split}
\end{equation}
Finally, substituting  (\ref{eq12}), (\ref{eq13}) into
(\ref{eq11}) gives the second-order difference equation for $L(n)$
\begin{equation*}
\begin{split}
4&(4n+3)(4n+5)(5n^2-2a^2-2b^2)L(n+1)+2(n+1)p(n)L(n) \\ &-n(n+1)
(5(n+1)^2-2a^2-2b^2)q(n)L(n-1)=0,
\end{split}
\end{equation*}
with initial conditions $L(0)=C(0),$
$$
L(1)=\left(\frac{1}{3}-\frac{2}{15}(a^2+b^2)\right)A(0)+
\left(\frac{1}{6}(a^2+b^2)-\frac{2}{15}(a^4+b^4
-4a^2b^2)-\frac{1}{30}\right)C(0),
$$
derived from (\ref{eq06}), (\ref{eq07}), (\ref{eq09}) and
polynomials $p(n), q(n)$ defined in (\ref{p}), (\ref{q}).
%\begin{equation*}
%\begin{split}
%p&(n)=30n^7+105n^6+n^5(145-52(a^2+b^2))+n^4(100-130(a^2+b^2))
%\\ &+n^3
%(35-124(a^2+b^2)+56(a^4+b^4)-208a^2b^2)+n^2(5-56(a^2+b^2)+84(a^4+b^4)
%\\ &-312a^2b^2)
%+n(80a^2b^2(a^2+b^2)-16(a^6+b^6)+48(a^4+b^4-3a^2b^2)-14(a^2+b^2))
%\\ &+(10(a^2-b^2)^2-2(a^2+b^2)+40a^2b^2(a^2+b^2)-8(a^6+b^6)),
%\end{split}
%\end{equation*}
%\begin{equation*}
%\begin{split}
%q(n)&=n^8-6n^6(a^2+b^2)+n^4(9(a^4+b^4)+30a^2b^2)-n^2(28a^2b^2(a^2+b^2)+4(a^6+b^6))
%\\ &+16a^2b^2(a^2-b^2)^2.
%\end{split}
%\end{equation*}

If we put $l(n)=L(n)/(n!)^4, n=0,1,2, \ldots,$ then it is easily
seen that the sequence $l(n)$ satisfies the following recurrence
equation:
\begin{equation*}
\begin{split}
4&(4n+3)(4n+5)(5n^2-2a^2-2b^2)n^3(n+1)^3l(n+1)+2n^3p(n)l(n) \\ &-
(5(n+1)^2-2a^2-2b^2)q(n)l(n-1)=0.
\end{split}
\end{equation*}
Its characteristic polynomial $64\lambda^2+12\lambda-1=0$ has two
different zeros $\lambda_1=-1/4,$ $\lambda_2=1/16,$ and by
Poincar\'e's theorem, we get
\begin{equation}
\lim_{n\to\infty}\left(\frac{|L(n)|}{(n!)^4}\right)^{\frac{1}{n}}\le
\frac{1}{4}. \label{lr}
\end{equation}
The limit inequality (\ref{lr}) implies
$$
\lim_{n\to\infty}\sum_{k=0}^{\infty}F(n,k)=0, \qquad
\lim_{k\to\infty}\sum_{n=0}^{\infty}G(n,k)=0
$$
and therefore, we have (see \cite[p.9]{moze}, \cite[\S 2]{ks})
$$
\sum_{k=0}^{\infty}F(0,k)=\sum_{n=0}^{\infty}G(n,0),
$$
yielding the theorem with $d_n=D(n-1),$ $L_n=L(n),$ $A_0=A(0),$
$B_0=B(0),$ $C_0=C(0).$     \qed

\section{Proof of Theorem \ref{t2}.}

To deduce (\ref{idt2}) from (\ref{eq04}), take $A(0)=C(0)=0,$
$B(0)=1$ and apply the following formula (see \cite[Corollary
2]{momo}):
\begin{equation}
\sum_{k=0}^{\infty}F(0,k)=\sum_{n=0}^{\infty}(F(n,n)+G(n,n+1)).
\label{last}
\end{equation}
Since in this case $L(n)=0$ for all $n\ge 0,$ an easy computation
of the right-hand side of (\ref{last}) by (\ref{eq06}),
(\ref{eq07}), (\ref{eq12})--(\ref{eq14}) and substitution
(\ref{ab})  lead to the desired conclusion.  \qed

\end{document}